\documentclass[runningheads]{llncs}
\usepackage[T1]{fontenc}
\usepackage{graphicx}
\usepackage{amsmath,amsfonts}
\usepackage{multirow}
\usepackage{authblk}
\usepackage{cite}
\usepackage{url}
\usepackage{pgf}
\usepackage{color}
\usepackage{todonotes}
\usepackage{booktabs} 
\usepackage{algorithm}
\usepackage{algpseudocode}
\usepackage{pgfplots}
\usepackage{caption}
\usepackage{soul}
\usepackage{bbm}
\usepackage{enumerate}
\usepackage{booktabs}

\usepackage{todonotes}

\begin{document}
\title{
Bound Tightening using Rolling-Horizon Decomposition for Neural Network Verification
}
\titlerunning{Bound Tightening using RH Decomposition for NN Verification}
%
\author{Haoruo Zhao\inst{1,2} \and
Hassan Hijazi\inst{2} \and
Haydn Jones\inst{2,3} \and Juston Moore\inst{2} \and Mathieu Tanneau\inst{1} \and Pascal Van Hentenryck\inst{1}}
\authorrunning{Zhao et al.}
%
\institute{Georgia Institute of Technology, Atlanta GA 30308, USA\\
\email{\{hzhao306,mathieu.tanneau\}@gatech.edu}, \email{pascal.vanhentenryck}@isye.gatech.edu\and
Los Alamos National Laboratory, Los Alamos NM 87545, USA\\
\email{\{hlh,hjones,jmoore01\}@lanl.gov}\and
University of Pennsylvania, Philadelphia PA 19104, USA\\
\email{haydnjonest@gmail.com}}
\maketitle              
\begin{abstract}
Neural network verification aims at providing formal guarantees on the output of trained neural networks, to ensure their robustness against adversarial examples and enable deployment in safety-critical applications.
This paper introduces a new approach to neural network verification using a novel mixed-integer programming (MIP) rolling-horizon decomposition method.
The algorithm leverages the layered structure of neural networks by employing optimization-based bound tightening (OBBT) on smaller sub-graphs of the original network in a rolling-horizon fashion and tightening the bounds in parallel.
This strategy strikes a balance between achieving tighter bounds and ensuring the tractability of the underlying mixed-integer programs.
Extensive numerical experiments, conducted on instances from the VNN-COMP benchmark library, demonstrate that the proposed approach yields significantly improved bounds compared to existing efficient bound propagation methods.
Notably, the proposed method proves effective in solving open verification problems.
Our code is built and released as part of the open-source mathematical modeling tool Gravity (\url{https://github.com/coin-or/Gravity}), which is extended to support generic neural network models.


\end{abstract}
\section{Introduction}

Neural networks are being applied in critical systems and high-consequence decision-making settings, e.g., power systems \cite{ARPAe} and autonomous driving \cite{bunel2018unified}.
How can we trust these models when the stakes are too high, when the price of failure is prohibitive or even life-threatening?
To justify trust, these models need to provide robustness guarantees. 
For example, in the context of power grid applications, a guarantee that a slight change in input (power system state) will not lead to unreasonable fluctuations in output (control actions predicted by the network). Mathematical optimization can provide such guarantees. While local optimization methods are acceptable for training these models, global optimality---i.e., a formal certificate that no better solution exists---is needed to provide robustness guarantees. There has been significant recent interest in providing verifiable properties for neural networks using global methods such as mixed-integer programming (MIP) \cite{gowaldual,dvijotham2020efficient,dathathri2020enabling}.
MIP solvers are appealing because they can, in principle, perform \emph{complete} verification for many network architectures.
Unfortunately, global optimization methods come with a hefty computational price tag, and using off-the-shelve solvers is not a scalable approach.
 
 Improving the scalability of global methods for neural network verification is an active research area, with approaches such as mathematical reformulations \cite{tsay2021partitionbased}, cutting planes \cite{anderson2019strong}, zonotopes \cite{kochdumper2023open}, custom-built branch-and-bound algorithms and relaxations \cite{zhang2018efficient,salman2019convex,xu2020automatic,xu2021fast,wang2021beta,zhang22babattack,zhang2022general,Zhang2023}, to name a few.
 In 2020, a community-driven effort led to the creation of the VNN competition \cite{brix2023first}, which has been held yearly since. The competition's goal is to ``allow researchers to compare their neural network verifiers on a wide set of benchmarks''.
The $\alpha,\beta$-CROWN team has consistently won this competition since its inception \cite{xu2020automatic,xu2021fast,wang2021beta,zhang22babattack,zhang2022general,Zhang2023}.
In this paper, we are hoping to bring the use of mixed-integer programming for neural network (NN) verification one step closer to viability. For this purpose, we tackle the verification problem as defined in \cite{brix2023first}, using a mixed-integer programming decomposition approach combined with optimization-based bound tightening (OBBT) \cite{caprara2010global}.
 OBBT is extensively used in global optimization solvers to reduce variable domains, especially for nonconvex mixed-integer nonlinear programs (MINLPs) \cite{gleixner2017three}. The OBBT algorithm solves two auxiliary optimization problems for each decision variable. In its initial form, OBBT relies on convexifying the feasible region and using variables' bounds as objective functions.
 In this work, we propose to preserve the mixed-integer nature of the subproblems, leveraging instead our proposed rolling horizon decomposition.
 One nice property of OBBT is its amenability to parallelization, since each auxiliary problem can be run independently. We take advantage of parallelization along with other speedup methods such as early termination using cutoff values as outlined in our approach below.

\section{Problem Statement}

In the domain of neural network verification, a ``white-box'' model is given, granting full visibility into the network's architecture and parameters \cite{brix2023first}.
Our verification challenge, based on the framework established by Bunel et al. \cite{bunel2018unified}, is to ascertain whether a neural network, denoted as function $f$ with $L$ layers, produces outputs that satisfy a desired property $P$ for all inputs within a specified range $\mathcal{C}$.

Formally, we verify that for any input $\mathbf{x}_{0} \in \mathcal{C}$, the network's output $\mathbf{y}^{(L)}$ adheres to the property $P(\mathbf{y}^{(L)})$, encapsulated by the implication:
$$
\mathbf{x}_{0} \in \mathcal{C} \Rightarrow P(\mathbf{y}^{(L)}).
$$

For instance, in assessing local robustness, we determine whether all inputs within an $\epsilon$-ball around a data point $a$ with label $y_a$ are classified as $y_a$ by the network.
This property is widely used in image classification cases, where we assess if a network's output label remains consistent under perturbations within a small tolerance.

\subsection{MIP Encoding for Trained Neural Networks}
Transitioning neural network architectures into a mixed-integer programming (MIP) format is key for verification. This process, in line with Tjeng et al.'s \cite{tjeng2017evaluating} methodology, enables the application of mathematical programming for thorough network analysis. The MIP model thus becomes a crucial tool for effective verification strategies.

\begin{figure}[H]
    \centering
    \includegraphics[width=\textwidth]{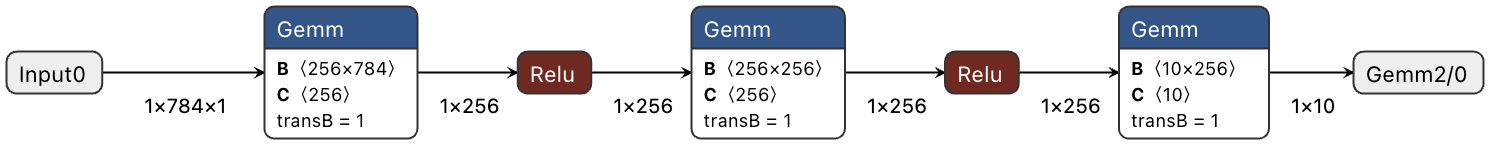}
    \caption[An Example of Fully Connected Neural Network with ReLU Activation]{A fully connected neural network with ReLU activation functions, showcasing the architecture used for MNIST digit classification.}
    \label{fig:mnist_fc}
\end{figure}

Consider a neural network with an input vector $ \mathbf{x}^{(0)} := \mathbf{x}_0 \in \mathbb{R}^{n_0} $. In this network, $ n_i $ represents the number of neurons in the $ i $-th layer. The network consists of $ L $ layers, where each layer $ i $ has an associated weight matrix $ W^{(i)} \in \mathbb{R}^{n_i \times n_{i-1}} $ and a bias vector $ b^{(i)} \in \mathbb{R}^{n_i} $, for $ i \in \{1, \ldots, L\} $. Let $ \mathbf{y}^{(i)} $ denote the pre-activation vector and $ \mathbf{x}^{(i)} $ the post-activation vector at layer $ i $, with $ \mathbf{x}^{(i)} = \sigma(\mathbf{y}^{(i)}) $. The output of the network is $ \mathbf{y}^{(L)} $. Although $ \sigma $ could be any activation function, we will assume it is the ReLU function throughout this paper. Figure~\ref{fig:mnist_fc} presents a fully connected neural network with ReLU activation functions.

In network verification, we typically define $f$ such that a non-negative outcome of $f(\mathbf{y}^{(L)}) \ge 0$ indicates the satisfaction of the desired property. Therefore, if an adversarial input $\mathbf{x} \in \mathcal{C}$ results in a negative value of $f$, this is taken as evidence that the verification instance fails to meet the required property, providing a counter-example. In short:
\[
\begin{cases}
    & \text{If $\exists$ } \mathbf{x_{\rm adv}^{(0)}} \in \mathcal{C} \text{ such that } f(\mathbf{y_{\rm adv}^{(L)}}) < 0, \text{then } P(\mathbf{y_{\rm adv}^{(L)}}) \text{ does not hold}.\\

    & \text{If $\forall$ } \mathbf{x^{(0)}} \in \mathcal{C}, f(\mathbf{y^{(L)}}) \geq 0, \text{then } P(\mathbf{y^{(L)}}) \text{ holds}.
\end{cases}
\]

The optimization problem is as follows:
\begin{equation}
\begin{aligned}
\min \quad & f(\mathbf{y}^{(L)}) \\
\text{s.t.} \quad & \mathbf{y}^{(i)} = W^{(i)}\mathbf{x}^{(i-1)} + b^{(i)}, & \forall i \in \{1, \ldots, L\}, \\
& \mathbf{x}^{(i)} = \sigma(\mathbf{y}^{(i)}), & \forall i \in \{1, \ldots, L-1\}, \\
& \mathbf{x}^{(0)} \in \mathcal{C}.
\end{aligned}
\label{eq:model1}
\end{equation}
Specifically, if all operators within the trained neural network are piecewise-linear, then the neural network can be linearly represented within the mixed-integer linear programming (MILP) framework. If the activation function is a ReLU function, the MILP formulation of $x = \sigma(y) = ReLU(y) = \max(0,y)$ is given by:
$$
x \geq 0, \ 
x \geq y, \ 
x \leq y - l(1 - z), \ 
x \leq u \cdot z, \ 
z \in \{0, 1\}
$$
where $l$, $u$ are the respective lower and upper bound on $y$. The binary variable $z$ indicates the activation state of the ReLU.

A ReLU neuron unit can be classified into different categories based on its input domain $[l, u]$: it is deemed ``Inactive'' if $u \leq 0$, ``Active'' if $l \geq 0$, and ``Unstabilized'' when $l < 0$ and $u > 0$. The neuron is called ``Stabilized'' if it meets either the Active or Inactive condition.

\subsection{Bound Tightening}
\label{sec:bound_tightening}
In practice, solving the mixed-integer program \eqref{eq:model1} can be computationally prohibitive for large neural networks. As a result, it has been proven effective to introduce additional bounds on intermediate layers. Specifically, the problem states
\begin{equation}\label{eq:extended_model}
\begin{aligned}
\min \quad & f(\mathbf{y}^{(L)}) \\
\text{s.t.} \quad& \text{Constraints of model }\eqref{eq:model1}, \\
\quad& \mathbf{x}^{(i)}_l \leq \mathbf{x}^{(i)} \leq \mathbf{x}^{(i)}_u, \quad \forall i \in \{1, \ldots, L-1\}, \\
\quad& \mathbf{y}^{(i)}_l \leq \mathbf{y}^{(i)} \leq \mathbf{y}^{(i)}_u, \quad \forall i \in \{1, \ldots, L\}.
\end{aligned}
\end{equation}
In model \eqref{eq:extended_model}, $\mathbf{y}^{(i)}_l$ and $\mathbf{y}^{(i)}_u$ denotes the vector lower and upper bound respectively for pre-activation output $\mathbf{y}^{(i)}$ at layer $i$; we similarly define these bounds for post-activation $\mathbf{x}^{(i)}$. Hence, the additional constraints form hyper-rectangles around the outputs.
It is crucial to note that the bounds on $\mathbf{y}^{(i)}$ play a much more important role than those on $\mathbf{x}^{(i)}$, as they serve as the input to the ReLU layer and consequently determine the stabilization of the ReLU neuron. When a ReLU is stabilized, the binary variable indicating its active state becomes fixed to either 0 (inactive) or 1 (active). This stabilization leads to a reduction in the number of binary variables in the problem.

There are several methods to derive bounds for intermediate layers:
\begin{itemize}
\item \textbf{Interval Bound Propagation (IBP):} This method employs interval bound propagation to establish bounds for each layer \cite{gowal2019effectiveness}.
\item \textbf{DeepPoly:} This method uses a custom polyhedral abstract domain relaxation. It assigns concrete lower and upper bounds to every neuron in a neural network. Symbolic bounds are formulated as linear combinations of the neurons in the network's previous layer.\cite{singh_2019}.
\item \textbf{CROWN:} This method efficiently leverages linear bound propagation to adaptively determine the lower and upper bounds of neural networks \cite{zhang2018efficient}.
\item \textbf{$\alpha$-CROWN:} This approach builds upon and enhances CROWN, further tightening the linear bounds by utilizing gradients \cite{xu2021fast}.
\item \textbf{Optimization-Based Bound Tightening (OBBT):} This approach represents the exact bound tightening technique, which involves solving MIP problems \cite{caprara2010global}.
\end{itemize}

The following remarks can be made about the above methods.
First, it is widely accepted that IBP, while simple and fast, yields very weak bounds, especially for deep networks.
Second, DeepPoly and CROWN are based on the same polyhedral relaxations, and therefore yield the same bounds \cite{salman2019convex}.
Third, because $\alpha$-CROWN augments CROWN with a gradient-based procedure, it achieves (both in theory and in practice) tighter bounds than CROWN \cite{salman2019convex,xu2021fast}.
Finally, as noted in \cite{salman2019convex}, because IBP, CROWN/DeepPoly and $\alpha$-CROWN are all based on polyhedral (linear) relaxations, they cannot break the convex relaxation barrier described in \cite{salman2019convex}, though $\alpha$-CROWN matches this theoretical limitation.
Note that the convex relaxation barrier presented in \cite{salman2019convex} is equivalent to LP-based OBBT, wherein bounds are tightened iteratively by solving only the linear relaxation of each OBBT problem.
Therefore, in order to further tighten bounds, any procedure must explicitly consider the binary variables associated to each ReLU neuron.

In our work, the primary objective is to obtain tighter bounds. To achieve this, we utilize optimization-based bound tightening, focusing on achieving the tightest box bounds for intermediate layers.
\section{Methodology}

We introduce the Optimization-Based Bound Tightening with Rolling Horizon (OBBT-RH) here. The method builds upon the MIP-based Optimization-Based Bound Tightening approach, as tighter bounds for intermediate layers in model \eqref{eq:extended_model} will in general result in faster solve times. Therefore, applying tighter bounds can speed up the verification process, allowing it to scale to larger models.

OBBT is formulated as two optimization subproblems for each neuron, seeking to find its maximum and minimum bound. Specifically, let $\mathbf{y}^{(t)}_k$ denote the $k$-th neuron at layer $t$ subject to network constraints. Given $0\leq s < t \leq L$, the problem states:
\begin{equation} \label{eq:obbt_instances}
\begin{aligned}
\text{max/min} \quad & \mathbf{y}^{(t)}_k \\
\text{s.t.} \quad & \mathbf{y}^{(i)} = W^{(i)}\mathbf{x}^{(i-1)} + b^{(i)}, & \forall i \in \{s+1, \ldots, t\}, \\
& \mathbf{y}^{(i)}_l \leq \mathbf{y}^{(i)} \leq \mathbf{y}^{(i)}_u, & \forall i \in \{s+1, \ldots, t\},\\
& \mathbf{x}^{(i)} = \sigma(\mathbf{y}^{(i)}), & \forall i \in \{s+1, \ldots, t-1\}, \\
& \mathbf{x}^{(i)}_l \leq \mathbf{x}^{(i)} \leq \mathbf{x}^{(i)}_u & \forall i \in \{s, \ldots, t-1\}.
\end{aligned}
\end{equation}
However, as the number of neural network layers included in the OBBT instances \eqref{eq:obbt_instances} increases, the OBBT process becomes intractable. To address this issue, we propose a decomposition method in OBBT-RH to reduce the number of layers considered in each OBBT instance.

The proposed OBBT-RH in Algorithm \ref{alg:OBBT-RH} leverages problem \eqref{eq:obbt_instances} to present an effective method for the verification of deep neural networks. OBBT-RH sequentially decomposes the neural network into manageable problem size, focusing on smaller sub-graphs of the original network. This method hence applies optimization-based bound tightening within a rolling horizon framework, effectively balancing the achievement of tighter bounds against maintaining the tractability of the mixed-integer programming problems involved.

\begin{algorithm}[!t]
\caption{OBBT-RH}
\label{alg:OBBT-RH}
\begin{algorithmic}[1]
\Require Rolling horizon sequence $\mathcal{S} = \{ (s_j, t_j) \mid j \in [J] \}$ where J is a positive integer
\For{$j = 1, \ldots, J$}
    \For{each neuron $k = 1, \ldots, n_i$}
        \State Solve problem \eqref{eq:obbt_instances} to obtain max and min values of $\mathbf{y}_k^{(t_j)}$.
    \EndFor
\EndFor
\Ensure Lower/upper bounds $\{\mathbf{y}^{(i)}_l\}, \{\mathbf{y}^{(i)}_u\}$ for all $i = 1, \ldots, L$.
\end{algorithmic}
\end{algorithm}

\subsection{Rolling Horizon Sequence}

\begin{figure}[!t]
    \centering
    \includegraphics[width=1.0\textwidth]{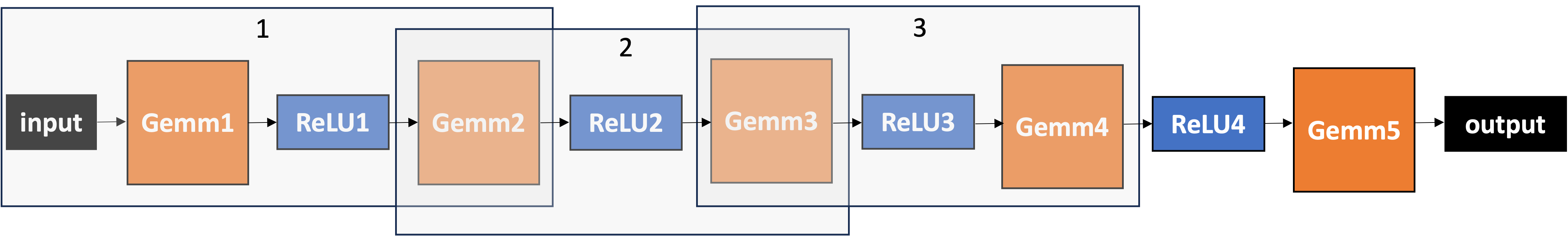}
    \caption{OBBT-RH with horizon length 2: $\mathcal{S} = \{(0,2), (1,3), (2,4)\}$, meaning that $(s_1, t_1) = (0, 2), (s_2, t_2) = (1, 3), (s_3, t_3) = (2, 4)$.}
    \label{fig:obbt_window}
\end{figure}

We introduce a rolling horizon strategy for selecting sequences $\mathcal{S}$ of layers within a neural network for bound tightening. The rolling horizon length, denoted as $ H $, determines the maximum number of General Matrix Multiplications (Gemm) layers to be included in each sub-graph. We begin with an initial pair $ (s_1, t_1) $ representing the starting and ending layers of the window, where $ s_1 $ is the starting layer and $ t_1 $ is the tightening layer. Each subsequent pair $(s_j, t_j)$ in the sequence $\mathcal{S}$ is constructed by shifting the window toward the end layer, with $H$ guiding the span of the window and indicating the layers selected for bound tightening in each iteration. For instance, considering the ReLU neural network, we selected $t_j = \min(j+1, L-1)$ and $s_j = \max(0, t_j-H)$. Indeed, the choice of $s_j$ and $t_j$ can be adjusted for various neural network architectures, highlighting the inherent flexibility of our approach.

The end goal is to tighten the neurons' bounds right before ReLU layers, thus each $ t $ in the sequence is a layer that precedes a ReLU layer. Figure \ref{fig:obbt_window} shows the OBBT-RH with horizon length $H = 2$ for a neural network with 5 layers. For $ H = 2 $, the sequence $\mathcal{S} = \{(0,2), (1,3), (2,4)\}$ includes pairs (\text{Input}, \text{Gemm2}), (\text{Gemm1}, \text{Gemm3}), and $(\text{Gemm2}, \text{Gemm4})$, each ending just before a ReLU layer. As $ H $ increases to 3, the sequence $\mathcal{S} = \{(0,2), (0,3), (1,4)\}$ includes pairs  $ (\text{Input}, \text{Gemm2}) $, $ (\text{Input}, \text{Gemm3}) $ and $ (\text{Gemm1}, \text{Gemm4}) $, and for $ H = 4 $, $\mathcal{S} = \{(0,2), (0,3), (0,4)\}$ and it includes $\{(\text{Input}, \text{Gemm2})$, $(\text{Input}, \text{Gemm3})$, $(\text{Input}, \text{Gemm4})\}$, now encompassing all Gemm layers leading to a ReLU operator. 

The length of the rolling horizon has direct impact on bound tightness and computational time:
A longer horizon considers more layers simultaneously, which could help in stabilizing more ReLU neurons and thus simplifying the MIP problem. However, this benefit comes at the cost of increased computational resources.
A shorter horizon will generally require less computation, beneficial for efficiency. However, it loses information in the previous layers and produces looser bounds.

One important advantage of this rolling horizon approach is the tightened bounds for neurons in previous layers can speed up the tightening process for neurons in subsequent layers. 
Note that IBP is equivalent to OBBT for the first ReLU layer if the intermediate layer is linear. We take advantage of this observation to avoid building sub-MIPs leading to this layer.

\subsection{Early Termination in the Case of ReLU Activation}
\paragraph{Maximizing Neuron Output:}
When solving the MIP \eqref{eq:obbt_instances} to maximize the output of a neuron $ y^{(t)}_k $, if the upper bound reaches zero, it indicates that the neuron will be inactive ($ y^{(t)}_k \leq 0 $). Hence, further bound tightening is unnecessary, and the process can be terminated early.

\paragraph{Minimizing Neuron Output:}
Similarly, when minimizing the output of a neuron, if the lower bound exceeds zero, it implies that the neuron will remain active ($ y^{(t)}_k > 0 $). In this case, the OBBT process can also be terminated early.

In either case, the ReLU neuron is stabilized. Therefore, introducing such an early termination criterion speeds up the bound tightening process.
In this paper, we simplify our model by assuming a linear structure for the neural network, wherein each layer is dependent solely on its immediate predecessor. However, real-world neural networks often exhibit more complex dependencies, where a given layer may depend not just on its immediate predecessor, but also on layers further back in the sequence. In such cases, the decomposition process becomes more complicated. We propose using Breadth-First Search (BFS) to navigate these complex dependencies. BFS can effectively trace the shortest path from the final layer $ t $ to an initial layer $ s $, ensuring that all relevant inter-layer dependencies, including those extending over multiple layers, are appropriately captured for conducting bound tightening of intermediate layers.

\subsection{Parallelization of OBBT sub-MIPs}
We take advantage of the graph-structure of neural networks by recognizing that the bounds on each neuron in a given layer can be independently tightened.
This observation allows us to distribute the computational effort of solving sub-MIPs corresponding to each neuron using independent threads in parallel.
Our implementation uses both the Message Passing Interface (MPI) for multi-machine clusters as well as single-machine multi-threading to parallelize our rolling-horizon algorithm.
\section{Numerical Results}

\begin{table}[!b]
\caption{Average number of neurons (inactive, active, stabilized, unstabilized) for each method in the benchmark.}
\centering
\begin{tabular}{@{}lcccc@{}}
\toprule
Method & Inactive & Active & Stabilized & Unstabilized\\ 
\midrule
IBP & 118.21 & 4.21 & 122.42 & 633.53 \\ 
CROWN & 435.67 & 18.72 & 454.39 & 301.56 \\ 
$\alpha$-CROWN & 466.67 & 19.93 & 486.60 & 269.35 \\
OBBT & 621.95 & 31.91 & 653.86 & 102.09 \\
OBBT-RH & \textbf{622.25} & \textbf{31.93} & \textbf{654.18} & \textbf{101.78} \\
\bottomrule
\end{tabular}
\label{table:neuron_states}
\end{table}
\subsection{Setup}

\paragraph{Benchmark dataset} Numerical experiments were conducted on all 90 instances from the \texttt{mnist\_fc} benchmark within the VNN-COMP benchmark library \cite{brix2023first}.

\paragraph{Comparison metrics} We primarily focus on the following metrics. First, the number of ReLU neurons stabilized (Table \ref{table:neuron_states}). Second, bounds range of ReLU neurons input (Table \ref{table:bounds_range}). Third, LP relaxation bounds of the MIPs with bounds on intermediate layers (Table \ref{table:lp_bounds}). Lastly, the proportion of instances that can be verified and their computing times (Table \ref{table:end2end}).

\paragraph{Implementation details} The \texttt{auto\_LiRPA} package is used to calculate IBP, CROWN, and $\alpha$-CROWN bounds. A GPU V100 is used for generating these bounds. Note that the code of \texttt{auto\_LiRPA} returns errors for five instances (51, 54, 60, 76, 85) when computing the $\alpha$-CROWN bounds. Therefore, for a fair comparison, these instances have been removed. 
The authors were not able to install DeepPoly and its required dependencies; this method is therefore not included in the results. 
Nevertheless, recall that DeepPoly and CROWN employ the same relaxations and thus yield numerically identical bounds \cite{salman2019convex}.
The optimization problems were formulated with mathematical modeling tool Gravity (\url{https://github.com/coin-or/Gravity}) in C++ and solved using Gurobi 10.0.2. The experiments are carried out on the High-Performance Computing (HPC) platform provided by the Partnership for an Advanced Computing Environment (PACE) at Georgia Institute of Technology, Atlanta, Georgia. Within this benchmark, the setting for OBBT-RH horizon length is dynamically adjusted based on the network architecture. Specifically, for neural networks with three layers, the rolling horizon length $H$ is set to two. In the case of networks with five layers, $H$ is set to three, while for those with seven layers, $H$ is set to five. OBBT can be seen as a special case of OBBT-RH, wherein the horizon length is equal to the depth of the neural network.

Additionally, during the initialization phase of OBBT-RH instances, IBP bounds are applied to intermediate layers. To further enhance efficiency, a parallel approach is employed for tightening the input of each ReLU layer. Given that each ReLU layer in this benchmark comprises 256 neurons, a total of 512 CPU threads are requested, with 2 threads used for each OBBT-RH instance. This strategy of parallelization enables the independent and simultaneous tightening of bounds for each ReLU layer, significantly boosting the process's overall efficiency. In addition, since the goal of each OBBT-RH instance is to tighten the bounds, we have set Gurobi's MIPFocus parameter to value 3 to enhance the quality of the bounds. Moreover, to avoid excessively long solving times for some instances, a 30-second time limit is enforced on each OBBT-RH instance. Therefore, with the need to solve both maximization and minimization problems for unstabilized neuron, the maximum time required for OBBT-RH in one layer is one minute. 

\subsection{Comparison with other bound tightening methods}

\begin{table}[!t]
\caption{Average bounds range for each method in the benchmark.}
\centering
\begin{tabular}{@{}lc@{}}
\toprule
Method & Bounds Range \\ 
\midrule
IBP & 2028.35 \\ 
CROWN & 96.16 \\ 
$\alpha$-CROWN & 53.52 \\ 
OBBT & 12.87 \\
OBBT-RH & \textbf{12.38} \\
\bottomrule
\end{tabular}
\label{table:bounds_range}
\end{table}

The aim is to compare the bound tightness produced by various methods. Specifically, the goal is to maximize the number of stabilized neurons, thereby reducing the count of binary variables. Additionally, for those unstabilized neurons, the approach seeks to tighten their input bounds as much as possible, aiding in the stabilization of subsequent layer neurons. After determining the bounds preceding the ReLU layers, the LP relaxation bounds of the MIPs will be compared with the intermediate layers' generated bounds. Lastly, the evaluation will focus on the end-to-end comparison of the number of verified instances as well as the time taken across different methods.

In terms of results, Table \ref{table:neuron_states} presents the average number of neurons categorized as inactive, active, stabilized and unstabilized for different bound tightening methods. The OBBT-RH method demonstrates a higher average across all layers from all instances.
Table \ref{table:bounds_range} shows the average bounds range for each verification method in a benchmark setting. OBBT-RH outperforms the other methods with the tightest bounds, indicating its effectiveness in bounding neuron activations within the network.
Table \ref{table:lp_bounds} presents the average LP relaxation bounds for MIPs using the bounds produced by each method. These LP relaxation bounds indicate the tightness of the solution space. OBBT shows the least negative value, implying a closer approximation to the MIP's optimal solution, whereas OBBT-RH has a slightly higher value than OBBT.
Table \ref{table:end2end} includes the end-to-end verification results by comparing the efficiency of MIP-based verification where variable bounds are computed using different bound tightening methods. BT stands for the bound tightening time and MIP stands for the final MIP solving time with the tightened bounds. A total time limit of 20 minutes is set for each approach. Overall, MIP with OBBT-RH bounds verifies the most instances (either by proving robustness or finding an adversarial example) than other approaches, and achieves a 2.2x speedup over CROWN-based methods. In addition, OBBT-RH yields 1.25x speedup over OBBT, mostly thanks to shorter bound-tightening times.

\begin{table}[!b]
\caption{Average LP relaxation bounds of the MIPs for each method in the benchmark.}
\centering
\begin{tabular}{@{}lc@{}}
\toprule
Method & LP Bounds \\ 
\midrule
IBP & -25232.25 \\ 
CROWN & -1169.45 \\ 
$\alpha$-CROWN & -577.91 \\
OBBT & \textbf{-0.64} \\ 
OBBT-RH & -24.77 \\
\bottomrule
\end{tabular}
\label{table:lp_bounds}
\end{table}

\begin{table}[!t]
\caption{End-to-end comparison on all instances. All time are in seconds.}
\centering
\begin{tabular}{@{}lccr@{\hspace{1em}}r@{\hspace{1em}}r}
\toprule
&&&\multicolumn{3}{c}{Time (sec)}\\
\cmidrule(lr){4-6}
Method & \#Verified & \#Timeout 
    & \multicolumn{1}{c}{BT} 
    & \multicolumn{1}{c}{MIP} 
    & \multicolumn{1}{c}{Total} \\
\midrule
IBP+MIP & 37 & 48 & \textbf{3.54} & 708.60 & 712.14 \\
CROWN+MIP & 62 & 23 & \textbf{3.54} & 345.89 & 349.43 \\
$\alpha$-CROWN+MIP & 61 & 24 & 4.93 & 340.74 & 345.66 \\
OBBT+MIP & 79 & 6 & 121.20 & 81.07 & 202.27 \\
OBBT-RH+MIP & \textbf{80} & \textbf{5} & 90.19 & \textbf{64.66} & \textbf{154.85} \\
\bottomrule
\end{tabular}
\label{table:end2end}
\end{table}

\begin{figure}[!t]
    \centering
    \includegraphics[width=1.0\textwidth]{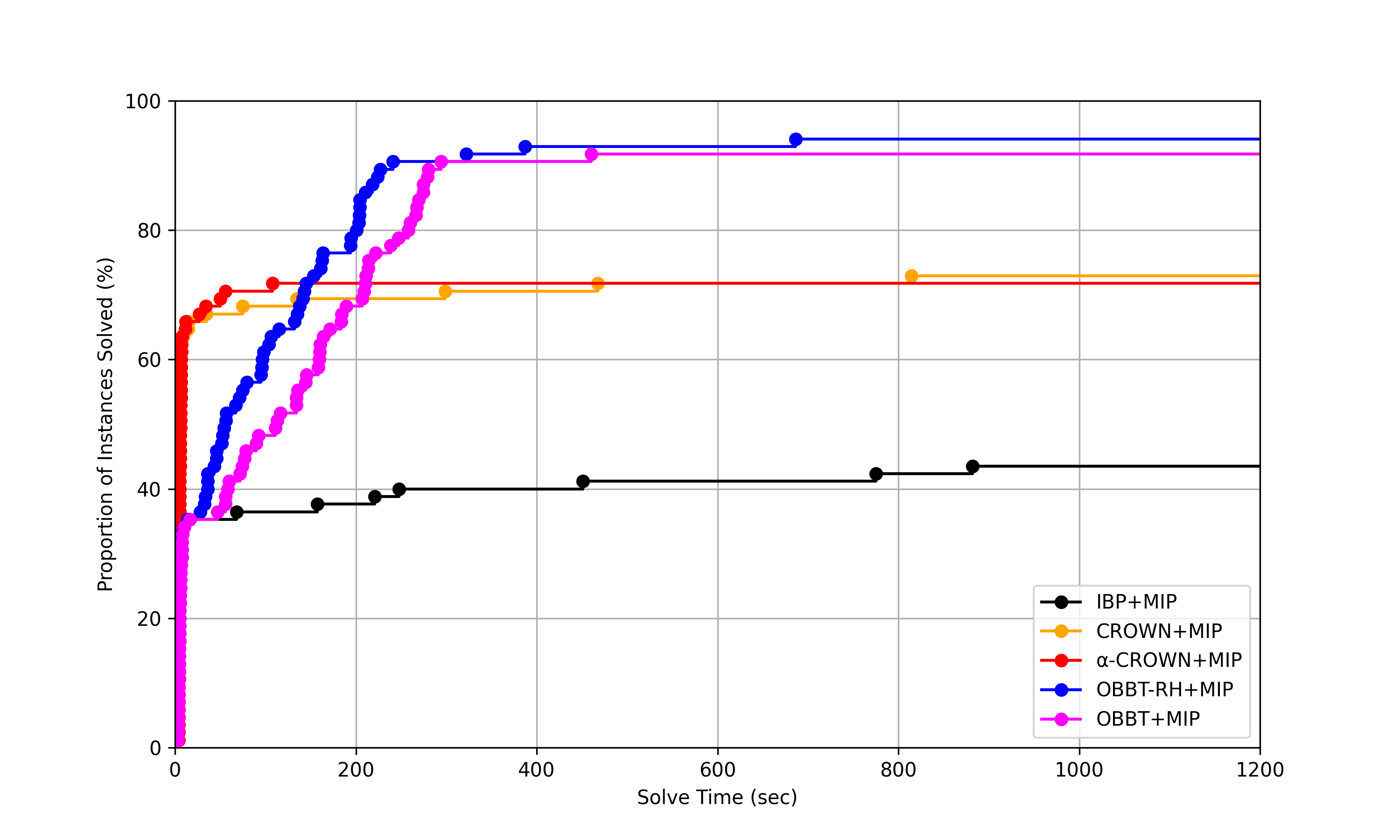}
    \caption{Performance comparison of the MIP-based verification for different bound-tightening methods.}
    \label{fig:cactus_plot}
\end{figure}

\subsection{Solving MIP with tight bounds}
After computing bounds for the intermediate layers, we use these bounds to solve the MIP instances. This is achieved by first obtaining the variable bounds for intermediate layers and then constructing the MIPs using the model as outlined in \eqref{eq:extended_model} for each verification instance. In the verification context, tightening the domain of the input for ReLU activations strengthens the MIP formulation by decreasing the big-M coefficients. Furthermore, it can stabilize ReLU units—eliminating binary variables—when a ReLU’s input is proven to be always non-negative or always non-positive. This combined effect—yielding stronger relaxations and reducing the number of binary variables—results in substantial improvements when solving the final MIPs corresponding to the verification problems with the tightened bounds.
In addition to incorporating the tight bounds into the MIP, we also set the Gurobi cutoff parameters to 0. This parameter setting is critical as it ensures that nodes with a lower bound greater than 0 are pruned, thereby enhancing the efficiency of the verification process. 

The end-to-end results are shown in Table \ref{table:end2end}, with OBBT-RH outperforming the other methods as discussed above.
Figure \ref{fig:cactus_plot} shows the performance of various methods in verifying instances within the benchmark, highlighting the trade-off between tight bounds and verification time. For low time budgets, fast bound tightening methods such as IBP, CROWN, and $\alpha$-CROWN are more effective, verifying a greater number of instances quickly. However, as the time budget increases, OBBT-RH and OBBT, which prioritize tighter bounds, begin to outperform the faster methods by verifying a higher proportion of instances. Notably, OBBT-RH is more efficient than OBBT, verifying more instances in less time.
Moreover, our approach has successfully closed several challenging instances (48, 78, 83) that previously posed difficulties for state-of-the-art complete verifiers. With our proposed method for generating OBBT-RH bounds, we successfully closed instance 48 within five minutes, a feat not achieved by any complete verifiers from the VNN competition.

\begin{figure}[!t]
    \centering
    \includegraphics[width=1.0\textwidth]{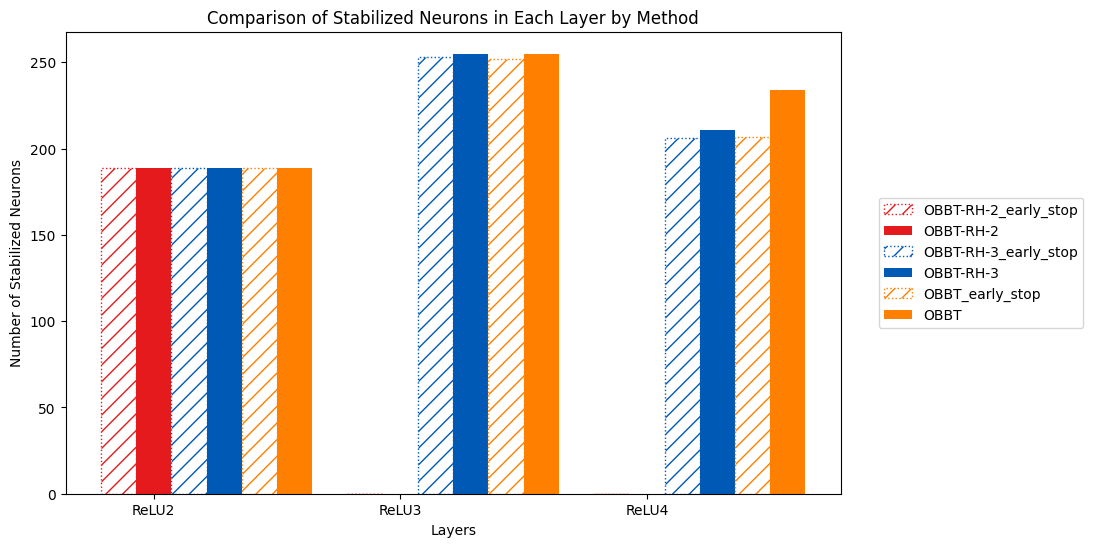}
    \caption{Comparison of stabilized neurons in each layer by different methods on instance 48.}
    \label{fig:stabilized_neurons}
\end{figure}

\subsection{Sensitivity Analysis 1: length of rolling horizon}

To assess the impact of horizon length on the performance of OBBT-RH, we conducted a sensitivity analysis. Due to the high computational cost of solving each OBBT instance to optimality, the sensitivity analysis was only conducted on instance 48. The analysis involved varying the horizon length of OBBT-RH, denoted as OBBT-RH-$H$, where $H$ represents the horizon length. For instance, OBBT-RH-3 indicates an OBBT-RH variant with a horizon length of 3. Figure \ref{fig:stabilized_neurons} compares the number of stabilized neurons across each layer by OBBT-RH with varying horizon lengths. The analysis reveals a clear trend: OBBT-RH-3 demonstrates a significant increase in the number of stabilized ReLU units compared to its counterpart, OBBT-RH-2, which employs a rolling horizon of length 2. OBBT-RH-2 stabilizes 0 ReLUs in the 3rd and 4th ReLU layers. It should also be emphasized that, considering early-stopping, OBBT-RH-3 stabilizes almost as many neurons as OBBT. This suggests that our proposed MIP-based approach, with an appropriately set rolling horizon, can effectively enhance the tightness of neuron bounds.

\begin{figure}[!t]
    \centering
    \includegraphics[width=1.0\textwidth]{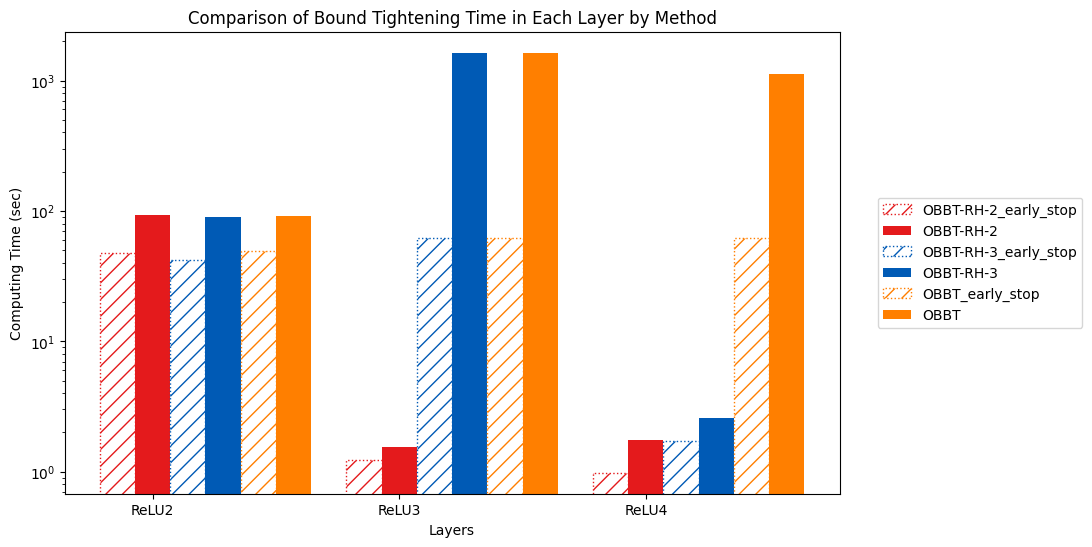}
    \caption{Comparison of bound tightening time for each layer by different OBBT-based methods on instance 48.}
    \label{fig:bound_tightening_time}
\end{figure}

\subsection{Sensitivity Analysis 2: the effect of early stop}
In the process of bound tightening, OBBT-RH is utilized to refine the input bounds of each neuron. This task is computationally intensive, requiring the solution of MIPs for both the upper and lower bounds of potentially hundreds of neurons.

Early stopping is a technique introduced to reduce the computational burden. It achieves this by terminating the bound tightening process before reaching the optimal solution under certain conditions.

To prevent any single OBBT-RH instance from taking an excessively long time, we have set a 30-second time limit. This decision is informed by empirical evidence indicating that many bounds exhibit only marginal improvement after 30 seconds. Therefore, this time limit is selected to balance the achievement of high-quality bounds with the need for reasonable computing time.

Figure \ref{fig:bound_tightening_time} illustrates the impact of early stopping on bound tightening time across different network layers, presented in a log-scale format. We cease further tightening of a neuron once it stabilizes, regardless of the application of early stopping. To assess the effectiveness of early stopping, we compare the number of stabilized neurons with and without this criterion. As depicted in the figure, early stopping—particularly with a 30-second time limit—does not significantly affect neuron stabilization. This finding supports the use of early stopping as an effective strategy to enhance the scalability of OBBT-RH without significantly compromising the quality of the bounds obtained. Due to the prolonged solving time required for tightening the inputs of the ReLU3 and ReLU4 layers with OBBT-RH-3 and OBBT without early stop (exceeding 4 hours), we resorted to using 24 threads to run a single OBBT instance. Therefore, the actual solving time with 2 threads would be even longer.
\section{Conclusion and Future Work}

We present a new algorithm integrating bound-tightening and mixed-integer programming with a rolling horizon strategy as a promising approach for verifying neural networks. Our method is particularly suited for handling networks where tighter bounds are needed but cannot be obtained from existing effective bound propagation methods. Future work will focus on extending this methodology to other types of nonlinear operators and further improving the computational efficiency of the verification process.

\section{Acknowledgements}
This research is partly funded by NSF award 2112533 and supported by the Laboratory Directed Research and Development program of Los Alamos National Laboratory under projects 20230578ER and 20240734DI.

\bibliographystyle{splncs04}
\bibliography{reference}
%






\end{document}